\documentclass[12pt]{article}

\usepackage{graphicx}
\usepackage{amsmath}
\usepackage{amsfonts}
\usepackage{amssymb}
\usepackage{ifthen}

\newtheorem{thm}{Theorem}[section]
\newtheorem{lemma}[thm]{Lemma}
\newtheorem{defn}[thm]{Definition}
\newtheorem{defns}[thm]{Definitions}

\newtheorem{example}[thm]{Example}
\newtheorem{examples}[thm]{Examples}

\newtheorem{propn}[thm]{Proposition}
\newtheorem{cor}[thm]{Corollary}
\newtheorem{prob}{Problem}

\newcommand{\os}{\le_O}
\newcommand{\ons}{\triangleleft_O}
\newcommand{\cs}{\le_C}
\newcommand{\cns}{\triangleleft_C}

\newcommand{\nons}{\not\ons}
\newcommand{\nos}{\not\os}

\newcommand{\gl}[2]{\mathop{GL}_{#1}(#2)}

\newcommand{\T}{\mathop{Tor}}

\newcommand{\Cr}{\prod}
\newcommand{\Dr}{\bigoplus}

\newcommand{\cen}{\mathop{C}}
\newcommand{\zg}{\mathop{Z}}

\newcommand{\norm}{\mathop{N}}
\newcommand{\Aut}{\mathop{Aut}}

\newcommand{\sub}{\mathcal{S}}

\newcommand{\hht}{\mathop{ht}}

\newcommand{\h}{\mathop{h}}

\newcommand{\df}[1]{\textbf{#1}}

\newenvironment{prf}[1][Proof]{\textbf{#1.}
}{\ \rule{0.5em}{0.5em} \vskip .5em}
\begin{document}

\title{Classifying Spaces of Subgroups \\ \hfill of Profinite Groups}



\date{November 2007}
\author{Paul Gartside and Michael Smith}

\maketitle

\begin{abstract}
The set of all closed subgroups of a profinite carries a natural
profinite topology. This space of subgroups can be classified up to
homeomorphism in many cases, and tight bounds placed on its complexity
as expressed by its scattered height. 
\end{abstract}

\section{Introduction}
In this paper we continue the task of topologically classifying the
space of closed subgroups of a profinite group. By
definition, a {\sl profinite} group $G$ is one that can be represented as
a {\sl projective} limit of {\sl finite} groups, $\lim_{\leftarrow}
G_{\lambda}$. Writing $\sub (G)$ 
for the set of closed subgroups of a profinite group,
$G=\lim_{\leftarrow} G_{\lambda}$, one sees that $\sub (G) =
\lim_{\leftarrow} \sub (G_\lambda)$. Giving the finite set
$\sub(G_\lambda)$ the discrete topology for each $\lambda$, we see
that the projective limit $\sub(G)$ picks up a natural topology. This
topology is profinite (compact, Hausdorff and zero-dimensional).
An alternative description of the topology on 
$\sub (G)$ is that it is the subspace topology inherited by $\sub(G)$ from the
space of all compact subsets of $G$ with the Vietoris topology, and so
the topology is independent of the particular projective
representation of $G$.

In previous work \cite{GS1} we have shown how to use the space of
subgroups to count the closed (and closed normal) subgroups of a
profinite group, and gave a partial classification (up to
homeomorphism) in certain cases. From that work we know to concentrate
on countably based (i.e. second countable) profinite groups. 

We know
that the space of subgroups $\sub (G)$ is countable if and only if $G$
is a central extension of $\Dr_{i=1}^k \mathbb{Z}_{p_i}$ where the
$p_i$s are distinct primes, and $\mathbb{Z}_{p_i}$ is the $p_i$-adic
integers. In this situation, $\sub (G) \cong \omega^k . n +1$. 

Recall that a profinite space is homeomorphic to the Cantor set if and
only if it is countably based and has no isolated points. Hence
$\sub(G)$ is homeomorphic to the Cantor space if and only if 
$\sub (G)$ has no isolated points; and we proved that $\sub (G)$ has
isolated points 
if and only if $G$ is finitely generated virtually nilpotent and only
finitely many primes divide the order of $G$. 

What remains is the case when we have
a countably based profinite group $G$, whose subgroup space is not
countable but does have isolated points. We tackle the problem of
classifying these subgroup spaces by analyzing the structure of the
isolated points via the Cantor-Bendixson process. 

Details of this
process, and relevant topological results, are given
in Section~\ref{bkgd}. In summary, for any space $X$ let $i(X)$ be the
set of isolated points 
in $X$, and call $X'=X \setminus i(X)$ the \df{derived set} of
$X$. We can then take the derived set of the derived set of $X$, and
so on (potentially transfinitely). For a countably based profinite
space this procedure will terminate at some least countable ordinal
$\alpha$, either with an empty derived set -- in which case the space
$X$ is countable -- or with a derived set homeomorphic to the Cantor set.
The \df{scattered height} of $X$ is the ordinal $\alpha$.

We show (Section~\ref{is_sol_sc}) that if $\sub (G)$ has any isolated
points then in fact all (the infinitely many) open subgroups are
isolated and they are dense in $\sub(G)$. Next we investigate the
properties of \df{solitary} subgroups -- those which are isolated in
$\sub (G)'$. From this we prove a very tight restriction on the
scattered height of a subgroup space (which is potentially any
one of the uncountable family of countable ordinals),
 namely that if $G$ is profinite
and $\sub (G)$ has any isolated points, then $\hht (\sub(G)) \le k+1$
where $k$ is the finite number of primes $p$ such that $p^\infty$
divides the order of $G$. In particular, if $G$ is pro-$p$ then
$\hht(\sub(G)) \le 2$. We show in what follows that all scattered
heights of $\sub (G)$ permitted by these results are attained.

Since a profinite group $G$ with isolated subgroups has an open normal
subgroup isomorphic to $N=\Dr_{i=1}^k G_{p_i}$ where each $G_{p_i}$ is a
pro-$p_i$ group, its space of subgroups has a clopen subspace
homeomorphic to $\sub (N) = \prod_{i=1}^k \sub (G_{p_i})$. Thus we
next focus on classifying subgroup spaces of pro-$p$
groups. 

Now suppose $G$ is a pro-$p$ group with isolated subgroups but an
uncountable subgroup space. A complete topological classification of
$\sub (G)$ seems well out of reach. Instead our results (in
Section~\ref{subg_prop}) show the
diversity of possibilities. 

Up to homeomorphism, there is a unique countably based profinite space
of scattered height $1$ with a dense set of isolated points. This
space is called the \df{Pe\l czy\'nski space}, $P$. We start, then, by
investigating when subgroup spaces are homeomorphic to $P$ -- in other
words those groups with no solitary subgroups. A wide variety
of pro-$p$ groups $G$ have $\sub (G) \cong P$ including: $\mathbb{Z}_p^n$
and the free pro-$p$ group on $n$ generators ($n>1$), insoluble
just-infinite pro-$p$ groups, and finitely generated nilpotent pro-$p$
groups with $h(G)>1$.

Next we consider the possibilities for $\sub (G)$ to have scattered
height $2$ -- in other words has solitary subgroups. One case is that
there are (countably) infinitely many  solitary subgroups. This occurs
if a non-abelian just-infinite pro-$p$ group is multiplied by
$\mathbb{Z}_p$. There are even poly-$\mathbb{Z}_p$ with $h(G)>1$ for
which this is true. We note that there are uncountably many pair-wise
non-homeomorphic countably based uncountable profinite spaces of
scattered height $2$. 

The second case is that there are just finitely many solitary
subgroups. Unlike the previous case, there is just a unique (up to
homeomorphism) countably based profinite space $X_n$ such that the
isolated points are dense, and there are exactly $n$ isolated points in $X_n'$.
This situation arises when the pro-$p$ group $G$ is virtually
$\mathbb{Z}_p$. In fact if $G$ is profinite with an open non-central
subgroup isomorphic to $\mathbb{Z}_p$, then $\sub(G) \cong P \oplus
(\omega . n+1)$ (so there are $n$ solitary subgroups). While if $G$ is
a non-abelian just-infinite pro-$p$ group then $\sub (G) \cong P \oplus
(\omega +1)$ (with exactly one solitary subgroup).

This leaves a natural open question: can an infinite finitely
generated pro-$p$ group which is not virtually $\mathbb{Z}_p$ have a
finite (non-zero) number of solitary subgroups? A natural conjecture
leading to a negative solution to this problem is refuted, and a
potential example exhibited, but not proved to have just finitely many
solitary subgroups.

We conclude (Section~\ref{gen_case}) with some
observations on the case of general profinite groups. In particular we
present examples showing that all permitted scattered heights are attained.

\section{Background Material}\label{bkgd}

\paragraph{The Cantor--Bendixon Process}\label{sctbdpr}
\begin{defn}
Let $X$ be a topological space. For $Y\subseteq X$, let $Y^{\prime}$ denote
the set of all limit points of $Y$, that is $Y\setminus Y^{\prime}$ is the
set of points of $Y$ which are isolated in $Y$. We define the following
transfinite sequence.

$X^{(0)}=X$,

$X^{(\alpha +1)}=(X^{(\alpha)})^{\prime}$, for $\alpha$ an ordinal,

$X^{(\lambda)}=\bigcap_{\mu <\lambda}X^{(\mu)}$, for $\lambda$ a limit ordinal.
\end{defn}
A space is \df{perfect} if it has no isolated points. 
\begin{lemma}\label{lctbdx}
Let $X$ be a Hausdorff space.
\begin{enumerate}
\item[(i)] $X^{(\alpha)}$ is closed in $X$ for every ordinal $\alpha$.
\item[(ii)] If $\alpha\leqslant\beta$ then $X^{(\alpha)}\supseteq X^{(\beta)}$.
\item[(iii)] $X^{(\alpha)}\setminus X^{(\alpha +1)}$ is the set of isolated
points of $X^{(\alpha)}$ for every ordinal $\alpha$, and is countable if
$X$ is countably based.
\item[(iv)] $(X^{(\alpha)})^{(\beta)}=X^{(\alpha +\beta)}$ for all ordinals
$\alpha$ and $\beta$.
\item[(v)] If $Y\subseteq Z\subseteq X$ then
$Y^{(\alpha)}\subseteq Z^{(\alpha)}$ for every ordinal $\alpha$.
\item[(vi)] There is a least ordinal $\lambda$ such that
$X^{(\lambda)}=X^{(\lambda +1)}$ and if $\alpha\geqslant\lambda$ then
$X^{(\alpha)}=X^{(\lambda)}$. $X^{(\lambda)}$ is perfect in itself. If $X$ is
countably based then $\lambda$ is countable.
\item[(vii)] If $X$ is compact, and $X^{(\alpha)}=\emptyset$ for some ordinal
$\alpha$, then there is a successor ordinal $\lambda$ such that
$X^{(\lambda)}=\emptyset$ and $X^{(\lambda -1)}$ is a finite non-empty
discrete space.
\item[(ix)] If $Y$ is open in $X$ then $Y^{(\alpha)}=Y\cap X^{(\alpha)}$ for
every ordinal $\alpha$.
\end{enumerate}
\end{lemma}
It follows from Lemma \ref{lctbdx} that every countably based Hausdorff
space can be written as the disjoint union of a countable set (the
\df{scattered part} of $X$) and a set which is perfect in itself (the
\df{perfect hull} of $X$); this is the Cantor--Bendixon theorem.
\begin{defns}
Let $X$ be a Hausdorff space. The \df{scattered height} of $X$, $\hht (X)$ is
the least ordinal $\lambda$ such that $X^{(\lambda)}=X^{(\lambda +1)}$.

Let $x\in X\setminus X^{(\hht (X))}$. $\hht (x,X)$ is defined to be the least
ordinal $\alpha$ such that $x\not\in X^{(\alpha +1)}$.
\end{defns}
Of course to say that $x\in X\setminus X^{(\hht (X))}$ has $\hht (X)=\alpha$
is precisely the same as saying that $x\in X^{(\alpha)}$ and that $x$ is
isolated in $X^{(\alpha)}$. The next lemma follows immediately from the
definitions.

\begin{lemma}\label{lhtopen}
Let $X$ be a Hausdorff space, $Y$ an open set in $X$ and $x\in Y$. Then
$\hht (Y)\leqslant \hht (X)$ and $\hht (x,X)=\hht (x,Y)$.
\end{lemma}

In relation  to Lemma \ref{lctbdx}(vii), a compact Hausdorff space $X$ for
which $X^{(\hht (X))}=\emptyset$ is sometimes known as a \df{scattered} space.
It is clear that a countably based profinite space is scattered if and only
if it is countable. On the other hand it is clear from Lemma \ref{lctbdx}
and Proposition \ref{pchcant} that if $X$ is an uncountable countably based
profinite space then $X^{(\hht (X))}$ is homeomorphic to the Cantor set. In
this sense, in order to understand the uncountable countably based profinite
spaces we need to understand how to `glue' a countable set onto a Cantor set.

\paragraph{Characterizations of Profinite Spaces}
We now describe some countably based profinite spaces and some well-known
characterisation theorems concerning them. These characterisation theorems
are all consequences of a celebrated theorem of Vaught; see, for example
Section 27 of \cite{Todor}.

Every ordinal is a topological space when considered with the order
topology. Successor ordinals, $\alpha +1$, are profinite, and
countably based if and only if $\alpha$ is countable. Note that
$\omega +1$ is homeomorphic to a convergent sequence.
\begin{lemma}
Let $\alpha$ be a non-zero countable ordinal, and $n$ be a positive integer.
Then $\omega^{\alpha}n+1$ is a countably infinite profinite space. Moreover
$(\omega^{\alpha}n+1)^{(\alpha)}=\{\omega^{\alpha},\omega^{\alpha}2,\ldots
,\omega^{\alpha}n\}$. In particular $\hht (\omega^{\alpha}n+1)=\alpha +1$
and $|(\omega^{\alpha}n+1)^{(\alpha)}|=n$.
\end{lemma}

\begin{propn}\label{pctprs}
Let $X$ be a countable profinite space. Then $X$ is homeomorphic to
$\omega^{\hht (X)-1}|X^{(\hht (X)-1)}|+1$.
\end{propn}

\begin{propn}\label{pchcant}
A space is homeomorphic to the Cantor set if and only if it a
countably based profinite space with no isolated points.
\end{propn}

\begin{defn}
Let $F_0=P_0=[0,1]$, the closed unit interval.

Let $F_1=F_0\setminus (1/3,2/3)$ and let $P_1=F_1\cup\{1/2\}$. That is to
form $F_1$ we have removed the open middle third interval of $F_0$, and
to form $P_1$ we have reinserted the mid-point of the deleted interval.

Continuing we remove the middle third intervals of the remaining segments
in $F_1$ to form $F_2$, that is $F_2=F_1\setminus ((1/9,2/9)\cup (7/9,8/9))$.
We let $P_2=F_2\cup\{1/6,5/6\}$.

We continue, forming $F_i$ by removing the middle third intervals
of the remaining segments in $F_{i-1}$ and forming $P_i$ by reinserting
the mid-points of the deleted intervals in $F_i$.

We let $F=\bigcap_{i=0}^{\infty}F_i$ and $P=\bigcap_{i=0}^{\infty}P_i$.
$P$ is called \df{Pe\l czy\'nski space}.
\end{defn}

Clearly $F$ is the standard construction of the middle third Cantor set,
 $P$ is a
profinite space and it is clear that $P\setminus F$ is the set of isolated
points in $P$ and that this is a countable dense subset of
$P$. Since$P^{\prime}=F$, $\hht (P)=1$. These facts characterise 
Pe\l czy\'nski space 
as a countably based profinite space.

\begin{propn}\label{pchpel}
Let $X$ be an infinite countably based profinite space. Then $X$ is
homeomorphic to Pe\l czy\'nski space, $P$, if and only if $\hht (X)=1$ and
$X\setminus X^{\prime}$ is dense in $X$.
\end{propn}

Pe\l czy\'nski space will arise frequently in the sequel as the space of
closed subgroups of a profinite group. Now Pe\l czy\'nski space has
scattered height $1$. We now briefly consider uncountable countably based
profinite spaces $X$ of scattered height $2$ and with
$X\setminus X^{\prime}$ dense in $X$. Now if
$X^{\prime}\setminus X^{(2)}$ is countably infinite then there are a lot
of possibilities for $X$ up to homeomorphism. But if
$X^{\prime}\setminus X^{(2)}$ is finite of size $n$ then $n$ determines
the space up to homeomorphism. 

\begin{propn}\label{pchpelconv}
Let $X$ be an uncountable countably based profinite space and let $n$
be a positive integer. Then $X$ is homeomorphic to
$P\oplus (\omega n+1)$ if and only if $\hht(X)=2$,
$X\setminus X^{\prime}$ is dense in $X$ and $|X^{\prime}\setminus X^{(2)}|=n$.
\end{propn}

\paragraph{Dimension, and Pro-$p$ Groups of Finite Rank}
In Section~\ref{subg_prop} we make frequent use of the dimension
function of pro-$p$ groups of finite rank. A profinite group $G$ has
\df{finite rank} if there is an integer $n$ such that every closed subgroup of
$G$ can be topologically generated by $\le n$ elements. It turns out
that pro-$p$ groups of finite rank are precisely the class of $p$-adic
analytic pro-$p$ groups. Hence $G$ has a \df{dimension}, $\dim(G)$, as a
$p$-adic analytic manifold. A more algebraic definition of $\dim(G)$ is
as the minimal size of a (topological) generating set of any powerful
open subgroup. A pro-$p$ group $G$ is \df{powerful} if $p$ is odd and
$G/\overline{G^p}$ is abelian or $p=2$ and $G/\overline{G^4}$ is
abelian.
The properties of $\dim$ of use here are as follows, see \cite{Dix}:
\begin{thm}\label{tdim}
Let $G$ be a pro-$p$ group of finite rank, $H \cs G$ and $K \cns
G$. Then (i) $\dim(H) \le \dim(G)$, (ii) $\dim (G) = \dim (K) + \dim
(G/K)$, (iii) $H$ is finite if and only if $\dim(H)=0$, and (iv) $H
\os G$ if and only if $\dim(H)=\dim(G)$.
\end{thm}

\paragraph{Subgroups of Products}
We frequently need to  analyse subgroups of direct products. For this
we use a profinite version 
of a well-known result for abstract groups. The abstract version has been
variously attributed to Goursat, Remak, and to Klein and Fricke. It is
clear from the proof of the result for abstract groups (see Theorem 1.61
of \cite{Schm} and Theorem 8.19 and p. 185 of \cite{Rose}), that the same
argument works for profinite groups.

\begin{propn}\label{pgour}
Let $G=G_1\times G_2$ be a profinite group and let $\pi$ and $\rho$ be the
corresponding projections onto $G_1$ and $G_2$ respectively. Let $H\cs G$.
Then $\pi (H)=HG_2\cap G_1$ and $H\cap G_1\cns HG_2\cap G_1$. Similarly
$\rho (H)=HG_1\cap G_2$ and $H\cap G_2\cns HG_1\cap G_2$. Also
$(H\cap G_1)\times (H\cap G_2)\cns H$. We have the following topological
isomorphisms:
$$ \frac{HG_2\cap G_1}{H\cap G_1}\cong\frac{HG_1\cap G_2}{H\cap G_2}\cong
\frac{H}{(H\cap G_1)\times (H\cap G_2)}.
$$
\end{propn}

\paragraph{The Space of Closed Subgroups} We can concretely describe
canonical basic open neighbourhoods of a subgroup 
in $\sub (G)$  for a profinite group $G$ as follows.
\begin{defn}
Let $G$ be a profinite group. For $H\cs G$ and $N\ons G$, we define
$B(H,N)=\{K\cs G\mid KN=HN\}$.
\end{defn}
\begin{lemma}\label{lclbase}
Let $G$ be a profinite group, and $H\cs G$. Suppose that
$(N_{\lambda})_{\lambda\in\Lambda}$ is a family of open normal subgroups of
$G$, forming a base for the open neighbourhoods of $1$ in $G$. Then
$(B(H,N_{\lambda}))_{\lambda\in\Lambda}$ forms a base for the open
neighbourhoods of $H$ in $\sub (G)$. 
\end{lemma}
\begin{cor}\label{csubwt}
Let $G$ be a profinite group. Then $\sub (G)$ is countably based if
and only if $G$ is countably based.
\end{cor}

A frequently used result, helping to reduce problems about the space
of subgroups of general
profinite groups to the space of subgroups of pro-$p$ groups, follows
from Proposition~\ref{pgour}.  
\begin{propn}\label{pctprfct}
If a profinite group $G$ is topologically isomorphic to $\Cr_{i=1}^n
G_{p_i}$ where $p_1, \ldots , p_n$ are distinct primes, and $G_{p_i}$
is a pro-$p_i$ group, then $\sub (G) \cong \prod_{i=1}^n \sub(G_{p_i})$.
\end{propn}

\section{Isolated and Solitary Subgroups, \\
\mbox{ }  \hfill Scattered  Height}\label{is_sol_sc} 

\paragraph{Isolated Subgroups} 

In \cite{GS1} we proved:
\begin{thm}[Gartside and Smith]\label{pchperf}
Let $G$ be a profinite group. Then the following are equivalent.
\begin{enumerate}
\item[(i)] $S(G)$ is not perfect.
\item[(ii)] $G$ has an isolated open subgroup.
\item[(iii)] Every open subgroup of $G$ is isolated.
\item[(iv)] The Frattini subgroup, $\Phi (G)$, the intersection of the
  maximal proper open subgroups, is an open normal subgroup of $G$.
\item[(v)] $G$ is finitely generated, virtually pronilpotent and only
finitely many primes divide the order of $G$.
\end{enumerate}
\end{thm}

From Proposition \ref{pchperf} we see that if $G$ is a profinite group with
$\sub(G)$ not perfect then the set of isolated subgroups of $G$ coincides
with the set of open subgroups of $G$. We make the following observation
about the set of open subgroups.

\begin{lemma}\label{losubsd}
Let $G$ be a profinite group. Then the set of open subgroups of $G$ is
dense in $\sub (G)$.
\end{lemma}

\begin{prf}
If $H\cs G$ and $N\ons G$ then $HN\os G$ and $HN\in B(H,N)$. Thus by Lemma
\ref{lclbase}, the set of open subgroups is dense in $\sub (G)$.
\end{prf}

We are thus led to the following description of $\sub (G)$ when $G$ is
a profinite group and $\sub (G)$ is not perfect.

\begin{propn}\label{pdsgnperf}
Let $G$ be a profinite group. If $\sub (G)$ is not perfect then the set
of open subgroups of $G$ is a countable open dense discrete subspace of
$\sub (G)$ and the set of closed non-open subgroups is a perfect set in
$\sub (G)$.
\end{propn}

\begin{prf}
If $G$ is a finitely generated
profinite group then $G$ has only countably many open subgroups. The result
now follows immediately from Theorem~\ref{pchperf} and Lemma~\ref{losubsd}.
\end{prf}

\begin{propn}
Let $G$ be a countably based profinite group. Then $\sub (G)$ is not perfect
if and only if the set of open subgroups of $G$ is open in $\sub (G)$.
\end{propn}

\begin{prf}
By Theorem~\ref{pchperf}, if $\sub (G)$ is not perfect then the set of
open subgroups of $G$ is open in $\sub (G)$.

Now suppose that the set of open subgroups of $G$, $O(G)$ is open in
$\sub (G)$. Then $O(G)$ is a Baire space. $O(G)=\bigcup_{H\os G}\{H\}$, and
each $\{H\}$ is closed in $O(G)$. Since $G$ is
countably based, $O(G)$ is countable and since $O(G)$ is a Baire space there
exists an open subgroup $H$ of $G$ such that $\{H\}$ is open in $O(G)$. But
then $\{H\}$ is open in $\sub (G)$; that is, $H$ is an isolated open subgroup
of $G$. Hence $\sub (G)$ is not perfect.
\end{prf}

\paragraph{Solitary subgroups}\label{ssolsg}

In this section we further analyse the space of closed subgroups of
profinite groups with isolated subgroups.

\begin{defn}
Let $G$ be a profinite group with $\sub (G)$ not perfect and $H\cs G$.
$H$ is said to be a \df{solitary subgroup} of $G$ if
$H\in\sub (G)^{\prime}$ and $H$ is isolated in $\sub (G)^{\prime}$.
\end{defn}

\begin{lemma}\label{ldisub}
Let $G$ be a profinite group with $\sub (G)$ not perfect, and let
$H\in\sub (G)^{\prime}$.
\begin{enumerate}
\item[(i)] $H$ is solitary in $G$ if and only if there exists an open normal
subgroup $N$ of $G$ such that if $K\cs G$ and $KN=HN$ then either $K=H$ or
$K\os G$.
\item[(ii)] If $K\os G$ with $H\leqslant K$ then $H$ is solitary in $G$ if and
only if $H$ is solitary in $K$.
\item[(iii)] If $K\cns G$ with $K\leqslant H$ and $H$ is solitary in $G$ then
$H/K$ is solitary in $G/K$.
\item[(iv)] Let $g\in G$. Then conjugation by $g$ induces a homeomorphism from
$\sub (G)^{\prime}$ to $\sub (G)^{\prime}$. In particular if $H$ is solitary in
$G$ then $H^g$ is solitary in $G$. Moreover if $N$ is an open normal subgroup
of $G$ such that  $B(H,N)\cap\sub (G)^{\prime}=\{H\}$ then
$B(H^g,N)\cap\sub (G)^{\prime}=\{H^g\}$.
\item[(v)] If $H$ is solitary in $G$ then $\Phi(H)\ons H$.
\item[(vi)] If $H$ is solitary in $G$ then there exists $K\os G$ with $H\cns K$
such that $K/H\cong\mathbb Z_p$ for some prime $p$.
\item[(vii)] $1$ is solitary in $G$ if and only if $G$ is virtually
$\mathbb Z_p$ for some prime $p$.
\end{enumerate}
\end{lemma}

\begin{prf}
By Lemma \ref{lclbase} $\{B(H,N)\cap\sub (G)^{\prime}\mid N\ons G\}$ is a base
for the open neighbourhoods of $H$ in $\sub (G)^{\prime}$. From this, (i)
follows immediately.

For (ii), let $K\os G$. Then $\sub(K)$ is open in $\sub(G)$.
Since $\sub (K)^{\prime}=\sub (K)\cap\sub
(G)^{\prime}$,  $\sub (K)^{\prime}$ is open in $\sub (G)^{\prime}$. 
The claim now follows from the fact that a point isolated in an open
subspace is isolated.

For (iii), let $K\cns G$ with $K\leqslant H$ and suppose that $H$ is solitary
in $G$. Then there exists an open neighbourhood $B$ of $H$ in $\sub (G)$ such
that $B\cap\sub (G)^{\prime}=\{H\}$. Let $\pi\colon G\to G/K$ be the natural
map. Then $\pi$ is an open continuous map and so by Proposition~5
 of \cite{FG2} $\sub (\pi)$ is an open continuous
map. Thus $B^*=\sub (\pi)(B)$ 
is open in $\sub (G/K)$. Clearly $H/K\in B^*\cap\sub (G/K)^{\prime}$ and so
$\{H/K\}\subseteq B^*\cap\sub (G/K)^{\prime}$. Let
$L/K\in B^*\cap\sub (G/K)^{\prime}$. Then $L\in B\cap\sub (G)^{\prime}=\{H\}$.
Thus $L=H$ and so $B^*\cap\sub (G/K)^{\prime}=\{H/K\}$ is open in
$\sub (G/K)^{\prime}$ as required.

For (iv), let $g\in G$. Then conjugation by $g$ is a topological automorphism
of $G$, and so conjugation by $g$ is a homeomorphism
from $\sub (G)$ to $\sub (G)$. Clearly if $K\cs G$ then
$K\in\sub (G)^{\prime}$ if and only if $K^g\in\sub (G)^{\prime}$. Also by
Proposition \ref{pdsgnperf} $\sub (G)^{\prime}$ is closed in $\sub (G)$.
Thus conjugation by $g$ is a homeomorphism from $\sub (G)^{\prime}$ to
$\sub (G)^{\prime}$. It is now clear that if $H$ is solitary in $G$ then $H^g$
is solitary in $G$. Now suppose that $N$ is an open normal subgroup of $G$ such
that $B(H,N)\cap\sub (G)^{\prime}=\{H\}$. Let
$L\in B(H^g,N)\cap\sub (G)^{\prime}$ and let $K=L^{g^{-1}}$. Then
$K\in\sub (G)^{\prime}$ and $K^gN=LN=H^gN$. So $KN=HN$ and hence $K=H$. Thus
$L=H^g$ as required.

For (v), suppose
that $\Phi(H)\nons H$. Let $N\ons G$. Then $N\cap H\os H$ and so
$N\cap H\not\leqslant\Phi(H)$. So there is a maximal open subgroup $M$ of
$H$ such that $N\cap H\nleqslant M$. Now $(N\cap H)M=H$ and so
$MN\cap H=H$. Thus $MN=HN$. But $M\in\sub (G)^{\prime}$ and $M\neq H$. So
by (i), $H$ is not solitary in $G$.

For (vi) suppose that $H$ is solitary in $G$. Then by (i) there is an
open normal subgroup $N$ of $G$ such that if $L\cs G$ and $LN=HN$ then
either $L=H$ or $L\os G$. Let $K=HN$. Clearly $K\os G$. If $k\in K$ then
$k=hn$ for some $h\in H$ and $n\in N$. Then $H^kN=H^{hn}N=H^nN=HN$. So
$H^k\in B(H,N)$. But $H^k\in S(G)^{\prime}$ and so $H^k=H$. Thus
$H\cns K$. Let $L/H\cs K/H$ with $L\neq H$. Clearly $L\in B(H,N)$. If
$L/H\nos K/H$ then $L\nos K$ and so $L\in\sub (G)^{\prime}$. But then
$L=H$, a contradiction. Hence $L/H\os K/H$. But now
$K/H\cong\mathbb Z_p$ for some prime $p$, because the $p$-adic
integers are characterised among all profinite groups as those whose
only non-open subgroup is the trivial subgroup.

Finally for (vii) suppose there is an open normal subgroup $N$ of $G$
with $N\cong\mathbb Z_p$ for some prime $p$. Then all subgroups of the
form $p^n\mathbb{Z}_p$ are isolated in $\sub (\mathbb{Z}_p)$, with $1$
as their unique limit. So $1$ is solitary in $\mathbb{Z}_p$ (and in
$N$). Then by (ii), $1$ is 
solitary in 
$G$. Conversely suppose that $1$ is solitary in $G$. Then by (vi) there
exists $K\os G$ with $K\cong\mathbb Z_p$ for some prime $p$. Clearly if
$N=K_G$ then $N\ons G$ and $N\cong\mathbb Z_p$.
\end{prf}

Note that by Lemma \ref{ldisub}(vi), if $G$ is an
infinite pro-$p$ group of finite rank and $H$ is a solitary subgroup of $G$
then $\dim (H)=\dim (G)-1$.

\paragraph{Bounding The Scattered Height}

\begin{lemma}\label{lprhtsub}
Let $G$ be a profinite group, $g\in G$, $H\cs G$ and $k$ be a non-negative
integer. If $\hht (H,\sub (G))=k$ then $\hht (H^g,\sub (G))=k$.
\end{lemma}

\begin{prf}
As in the proof of Lemma \ref{ldisub}(iv) conjugation by $g$ is a
homeomorphism from $\sub (G)$ to $\sub (G)$. We now proceed by induction
on $k$. Clearly if $k=0$ then $H\os G$ and so $H^g\os G$ and $\hht (H^g)=0$.
So suppose the result holds for all closed subgroups $K$ with
$\hht (K,\sub (G))=k$. Then if $K\cs G$ then $K\in\sub (G)^{(k+1)}$ if and
only if $K^g\in\sub (G)^{(k+1)}$. Also by Lemma \ref{lctbdx},
$\sub (G)^{(k+1)}$ is closed in $\sub (G)$. Thus conjugation by $g$ is a
homeomorphism from $\sub (G)^{(k+1)}$ to $\sub (G)^{(k+1)}$. Now let
$H\cs G$ with $\hht (H,\sub (G))=k+1$. Then $H\in\sub (G)^{(k+1)}$ and $H$
is isolated in $\sub (G)^{(k+1)}$. So clearly using the above
homeomorphism, $H^g$ is isolated in $\sub (G)^{(k+1)}$. That is
$\hht (H,\sub (G))=k+1$ as required.
\end{prf}

\begin{thm}\label{bd_sc_ht}
Let $G$ be a profinite group with $\sub (G)$ not perfect. Let $k$ equal the
number of primes $p$ such that $p^{\infty}\mid o(G)$. Then
$ht(\sub (G))\leqslant k+1$. In particular if $G$ is any pro-$p$ group then
$ht(\sub (G))\leqslant 2$.
\end{thm}

\begin{prf}
Note that by Theorem~\ref{pchperf}, $n$ is finite. Suppose for a
contradiction that there is a profinite group $G$ with $\sub (G)$ not
perfect and $\hht (\sub (G))>k+1$. Then there exists $H\cs G$ with
$\hht (H,\sub (G)=k+1$. So $H\in\sub (G)^{(k+1)}$ and $H$ is isolated in
$\sub (G)^{(k+1)}$. Hence by Lemma \ref{lclbase} there exists an open normal
subgroup $N$ of $G$ such that $B(H,N)\cap\sub (G)^{(k+1)}=\{H\}$. Let $K=HN$.
Then as in the proof of Lemma \ref{ldisub}(vi) if $k'\in K$ then
$H^{k'}\in B(H,N)$. But by Lemma \ref{lprhtsub} $H^{k'}\in\sub (G)^{(k+1)}$. So
$H^{k'}=H$. Hence $H\cns K$.

We now consider $\sub (K/H)$. Now $H$ is isolated in $\sub (K)^{(k+1)}$ and
so the trivial subgroup of $K/H$ is isolated in $\sub (K/H)^{(k+1)}$, that
is $\hht (1,\sub (K,H))=k+1$. Let $L/H\cs K/H$ with $L\neq H$. Then
$\hht (L/H,\sub (K,H))\leqslant k$. Thus $\sub (K/H)$ is countable and
$\hht (\sub (K/H)=k+2$. But by hypothesis the number of primes $p$ such that
$p^{\infty}\mid o(K/H)$ is at most $k$. So by
 of \cite{GS1},
$\hht (\sub (K/H))\leqslant k+1$, a contradiction.
\end{prf}

\section{Subgroup Spaces of pro-$p$ Groups}\label{subg_prop}

\paragraph{When is $\sub (G)$ Homeomorphic to Pe\l czy\'nski's Space?}
Lemma \ref{ldisub} provides us with the techniques to show that many classes
of pro-$p$ groups $G$ have $\sub (G)$ homeomorphic to Pe\l czy\'nski space.
Of course we have a useful characterisation of such groups using Proposition
\ref{pchperf}, Proposition \ref{pdsgnperf} and Proposition \ref{pchpel}.

\begin{propn}\label{pproppel}
Let $G$ be a pro-$p$ group. Then $\sub (G)$ is homeomorphic to Pe\l czy\'nski
space if and only if $G$ is infinite, finitely generated and has no solitary
subgroups.
\end{propn}

\begin{examples}\label{ezpnfree}
Let $p$ be a prime and let $n$ be an integer $>1$.
\begin{enumerate}
\item[(i)] $\sub (\mathbb Z_p^n)$ is homeomorphic to Pe\l czy\'nski space.
\item[(ii)] If $G$ is a free pro-$p$ group of rank $n$ then $\sub (G)$ is
homeomorphic to Pe\l czy\'nski space.
\end{enumerate}
\end{examples}

\begin{prf}
For (i), firstly note that $1$ is not solitary in $\mathbb Z_p^n$ by Lemma
\ref{ldisub}(vii). Now let $H\in\sub (\mathbb Z_p^n)^{\prime}$ with $H\neq 1$
and let $N\os\mathbb Z_p^n$. Since $H$ is not isolated, it is not
open, so $\dim (H) < \dim(\mathbb{Z}_p^n)$, and there is a positive
integer $m<n$ such that $H$ can be generated (topologically) by
$m$ elements. Also $HN\neq H$. Hence there exist $h_1,\ldots ,h_m\in H$ and
$n_1,\ldots ,n_m\in N$ with $n_1\not\in H$ such that
$\overline{\langle h_1n_1,\ldots ,h_mn_m\rangle}N/N=HN/N$. Let
$K=\overline{\langle h_1n_1,\ldots ,h_mn_m\rangle}$. Now $KN=HN$, and
as $\dim(K) < \dim(\mathbb{Z}_p^n)$, $K$ is not open, so not isolated,
and $K\in\sub (\mathbb Z_p^n)^{\prime}$. But
$K\neq H$. So by Lemma \ref{ldisub}(i), $H$ is not solitary in $\mathbb Z_p^n$.
Hence by Proposition \ref{pproppel}, $\sub (\mathbb Z_p^n)$ is homeomorphic
to Pe\l czy\'nski space.

Now for (ii) let $G$ be a free pro-$p$ group of rank $n$. For the definition
and basic properties of free pro-$p$ groups see Chapter 5 of \cite{WilPrf}.
Again $1$ is not solitary in $G$, and so let $H\in\sub (G)^{\prime}$ with
$H\neq 1$. As before by Proposition \ref{pproppel} it suffices to show that
$H$ is not solitary in $G$. So suppose for a contradiction that $H$ is
solitary in $G$. Then by Lemma \ref{ldisub}(vi) there exists $K\os G$ with
$H\cns K$ and $K/H\cong\mathbb Z_p$. Also by Lemma \ref{ldisub}(ii) $H$ is
solitary in $K$. By Theorem 5.4.4 of \cite{WilPrf} $K$ is a free pro-$p$
group of rank $m$ where $m$ is an integer $>1$. Since $K/H$ is abelian,
$H\geqslant K^{\prime}$. By Proposition 5.1.5 and Proposition 5.1.3$^{\prime}$
of \cite{WilPrf}, $K/K^{\prime}\cong\mathbb Z_p^m$. But by Lemma
\ref{ldisub}(iii) $H/K^{\prime}$ is solitary in $K/K^{\prime}$. This
contradicts (i). Thus $H$ is not solitary in $G$ as required.
\end{prf}

\begin{propn}\label{ppicore}
Let $G$ be an infinite finitely generated pro-$p$ group, and $H$ be a solitary
subgroup of $G$. Then $G/H_G$ has an open normal subgroup, topologically
isomorphic to $\mathbb Z_p^n$ for some positive integer $n$ where $n$ is at
most the (necessarily finite) number of conjugates of $H$ in $G$.
\end{propn}

\begin{prf}
By Lemma \ref{ldisub}(vi) $\norm_G(H)\os G$ and so $H$ has only finitely
many conjugates in $G$. Let $H_1,\ldots ,H_k$ be the conjugates of $H$ in
$G$. Now by Lemma \ref{ldisub}(iv) and the proof of Lemma \ref{ldisub}(vi)
there exists $N\ons G$ such that for every $i$, $H_i\cns H_iN$ and
$H_iN/H_i\cong\mathbb Z_p$. Clearly $H_G\leqslant (HN)_G$, $(HN)_G\ons G$
and for every $i$, $(H_iN)_G=(HN)_G$. Also for every $i$,
$(HN)_G/(H_i\cap (HN)_G)\cong ((HN)_GH_i)/H_i$ which is a non-trivial
closed subgroup of $H_iN/H_i$. So for every $i$,
$(HN)_G/(H_i\cap (HN)_G)\cong\mathbb Z_p$. Now
$$
\frac{(HN)_G}{H_G}=\frac{(HN)_G}{\bigcap_{i=1}^k(H_i\cap (HN)_G)}
\hbox{ which embeds in }
\Dr_{i=1}^k\frac{(HN)_G}{H_i\cap (HN)_G}\cong\mathbb Z_p^k.
$$
Thus $(HN)_G/H_G\cong\mathbb Z_p^n$ for some positive integer $n$ with
$n\leqslant k$, as required.
\end{prf}

\begin{cor}
Let $G$ be an insoluble just-infinite pro-$p$ group. Then $\sub (G)$ is
homeomorphic to Pe\l czy\'nski space.
\end{cor}

\begin{prf}
Suppose for a contradiction that $\sub (G)$ is not homeomorphic to
Pe\l czy\'nski space. Then by Proposition \ref{pproppel}, $G$ has a
solitary subgroup $H$. Now as $G$ is just-infinite and $H\nos G$,
$H_G=1$. But by Proposition \ref{ppicore}, $G$ is then soluble, a
contradiction.
\end{prf}

A profinite group $G$ is \df{poly--procyclic} if it has a finite
series of closed subgroups $1 = G_0 \cs G_1 \cs \cdots \cs G_n=G$ such
that $G_{i-1} \cns G_i$ and $G_i/G_{i-1}$ is procyclic (topologically
generated by one element).

Let $G$ be a poly--procyclic group. For each prime $p$, the
\df{$\mathbb{Z}_p$-length} of $G$, denoted $h_p(G)$,  is the number of
factors having Sylow $p$-subgroups isomorphic to $\mathbb{Z}_p$ in a
series in $G$ with procyclic factors. If $G$ is pro-$p$ then we write
$h(G)$ for $h_p(G)$, and refer to the $\mathbb{Z}_p$-length as the
\df{Hirsch length}. It turns out that for a poly--procyclic pro-$p$
group $G$, $h(G)=dim(G)$. 
\begin{propn}\label{psfgnp}
Let $G$ be a finitely generated nilpotent pro-$p$ group with $\h (G)>1$.
Then $\sub (G)$ is homeomorphic to Pe\l czy\'nski space.
\end{propn}

\begin{prf}
We first show that without loss of generality we may assume that $G$ is
torsion-free. Suppose for a contradiction that $\sub (G)$ is not homeomorphic
to Pe\l czy\'nski space but that the result holds for finitely generated
torsion-free nilpotent pro-$p$ groups. Again by Proposition \ref{pproppel},
$G$ has a solitary subgroup $H$. By Lemma \ref{ldisub}(vi) there exists
$K\os G$ with $H\cns K$ and $K/H\cong\mathbb Z_p$. Let $T=\T
(H)$. Then $T$ is finite. (This is the pro-$p$ analogue of the well
known fact that every finitely generated nilpotent abstract group has
a finite torsion subgroup. The proof is by induction on the Hirsch
length). 
Also $T\cns K$. Clearly $K/T$ is a finitely
generated torsion-free nilpotent pro-$p$ group. Since $K$ is open in
$G$ and $T$ finite, $\h (K/T)=dim (K/T) = \dim (K) = \dim (G)=\h
(G)>1$. 
But by Lemma \ref{ldisub}(iii) $H/T$ is solitary in
$K/T$, a contradiction.

So we now assume that $G$ is torsion-free. We proceed by induction on
$\h (G)$. Again assume for a contradiction that $G$ has a solitary subgroup
$H$. Then, again there exists $K\os G$ with $H\cns K$ and
$K/H\cong\mathbb Z_p$. Since $K$ is torsion-free and nilpotent
$H\zg (K)/\zg (K)$ is torsion-free (see 5.2.19 of \cite {Rob} for example).
Thus $H/(H\cap \zg (K))\cong H\zg (K)/\zg (K)$ is torsion-free.

Clearly $0\leqslant\h (H/(H\cap\zg (K)))\leqslant\h (G)-1$. We now exclude
the end cases. If $\h (H/(H\cap\zg (K)))=0$ then $H/(H \cap \zg (K))$ is
finite, and
$H\cap \zg (K)\os H$, so $H$ is central in $K$. Thus $K=H\times K/H$. But
$H$ is abelian and so $K\cong\mathbb Z_p^{\h (G)}$. By Example
\ref{ezpnfree}(i) $H$ is not solitary in $K$ and so by Lemma \ref{ldisub}(ii)
$H$ is not solitary in $G$, a contradiction. Hence $\h (H/(H\cap \zg (K)))>0$.
If $\h (H/(H\cap \zg (K)))=\h (G)-1$ then $\h (H\cap \zg (K))=0$ and so
$H\cap \zg (K)=1$, a contradiction as $K$ is nilpotent (see 5.2.1 of
\cite{Rob} for example). So $\h (H/(H\cap \zg (K)))\leqslant\h (G)-2$, and
$\h (G)>2$.

Hence there exists $L\cs K$ with $H\cap \zg (K)\leqslant L$ such that
$K/(H\cap\zg (K))=H/(H\cap\zg (K))\rtimes L/(H\cap\zg (K))$ and
$L/(H\cap\zg (K))\cong K/H\cong\mathbb Z_p$. Thus
$1<\h (K/(H\cap\zg (K)))<\h (G)$ and $K/(H\cap\zg (K))$ is a finitely
generated torsion-free nilpotent pro-$p$ group. Hence by the inductive
hypothesis $\sub (K/(H\cap \zg (K))$ is homeomorphic to Pe\l czy\'nski
space. But by Lemma \ref{ldisub}(ii,iii) $H/(H\cap\zg (K))$ is solitary in
$K/(H\cap \zg (K))$, a contradiction. Hence $H$ is not solitary in $G$, a
contradiction.
\end{prf}

In this section we have been attempting to solve:
\begin{prob}\label{pbchpel}
Characterise algebraically, the pro-$p$ groups $G$ with $\sub (G)$
homeomorphic to Pe\l czy\'nski space.
\end{prob}
From the above, we know the class
of groups we hope to characterise here is very diverse; it includes
uncomplicated groups like $\mathbb Z_p^2$, and also very complicated groups
like the Nottingham group.

Perhaps this question is more tractable for pro-$p$ groups of finite
rank. In particular we ask:
\begin{prob}
If $G$ is a uniform pro-$p$ group with $\dim (G)>1$, is $\sub (G)$
homeomorphic to Pe\l czy\'nski space?
\end{prob}
A positive answer to this question seems reasonable, given that the group
operation in a uniform pro-$p$ group of dimension $n$ can be `smoothed
out', to give a new group which is topologically isomorphic to
$\mathbb Z_p^n$ (see  \S 4.3 of
\cite{Dix}). Also a positive answer to this question would imply that if
$G$ is a pro-$p$ of finite rank and dimension $>1$ then $G$ has an open
normal subgroup $N$ such that $\sub (N)$ is homeomorphic to
Pe\l czy\'nski space.

\paragraph{An Almost Complete Classification for Nilpotent pro-$p$ Groups}
We now have done enough to describe $\sub (G)$ for $G$ any nilpotent pro-$p$
group with $w(G)\leqslant\aleph_1$ where $w(G)$, the \df{weight} of
$G$, is the cardinality of a minimal sized base of $G$. 
We can say which spaces can occur and
characterise when they occur. First  an easy lemma.

\begin{lemma}\label{lnilctvzp}
Let $G$ be a pronilpotent group which is virtually $\mathbb Z_p$ for some
prime $p$. Then $\sub (G)$ is countably infinite if and only if $G$ is
nilpotent.
\end{lemma}

\begin{prf}
Let $N\ons G$ with $N\cong\mathbb Z_p$. Suppose that $\sub (G)$ is
countably infinite. Then by Theorem~3.7 of \cite{GS1} we know that $N$
is central in $G$. But $G/\zg (G)$ is then nilpotent and so $G$ is
nilpotent. Conversely suppose that $G$ is nilpotent. Then $N\cap\zg (G)$
is non-trivial (see 5.2.1 of \cite{Rob} for example). So
$N\cap\zg (G)\os N$ and thus $N\cap\zg (G)\cong\mathbb Z_p$. Hence $G$ has
a central open subgroup topologically isomorphic to $\mathbb Z_p$. Thus
$\sub (G)$ is countably infinite by Theorem~3.7 of \cite{GS1}.
\end{prf}

\begin{thm}\label{tchnlpprp}
Let $G$ be a nilpotent pro-$p$ group with $w(G)\leqslant\aleph_1$. Then
$\sub (G)$ lies in precisely one of the following four classes. We give
specific characterisations of when $\sub (G)$ is in each class in terms
of properties of $G$.
\begin{enumerate}
\item[(i)] $\sub (G)$ is a finite discrete space precisely when $G$ is finite.
\item[(ii)] $\sub (G)$ is homeomorphic to $\omega n+1$ for some positive
integer $n$ precisely when $G$ is virtually $\mathbb Z_p$.
\item[(iii)] $\sub (G)$ is homeomorphic to Pe\l czy\'nski space precisely
when $G$ is finitely generated and $\h (G)>1$.
\item[(iv)] $\sub (G)$ is homeomorphic to $\{0,1\}^{w(G)}$ precisely when
$G$ is not finitely generated.
\end{enumerate}
\end{thm}

\begin{prf}
The result follows immediately from  Lemma
\ref{lnilctvzp}, Proposition \ref{pchperf}, and Proposition
\ref{psfgnp}, and
Corollary~4.3
, Corollary~5.7
 and Theorem~6.8
from our earlier paper on subgroup spaces of
profinite groups, \cite{GS1}.
\end{prf}

\paragraph{Scattered Height $=2$, Infinitely Many Solitary Subgroups}

We now show that there are pro-$p$ groups with subgroup spaces which
have scattered height $2$
and infinitely many isolated points in the derived set.
There are even such
examples amongst poly-procyclic pro-$p$ groups. In order to see this we need
to analyse subgroups of direct products.

\begin{lemma}\label{lsdpzp}
Let $H$ be a finitely generated pro-$p$ group, $Z\cong\mathbb Z_p$ and
$G=H\times Z$. Then $H$ is solitary in $G$ if and only if
$H^{\prime}\os H$.
\end{lemma}

\begin{prf}
Firstly suppose that $H$ is solitary in $G$. Then by Lemma
\ref{ldisub}(iii) $H/H^{\prime}$ is solitary in $G/H^{\prime}$. Now
$G/H^{\prime}$ is a finitely generated abelian pro-$p$ group. If
$H^{\prime}\nos H$ then $\h (G/H^{\prime})>1$ contradicting Proposition
\ref{psfgnp}. Thus $H^{\prime}\os H$. It is not hard to give a direct
argument avoiding the use of Proposition \ref{psfgnp} by considering
the torsion subgroup of $H/H^{\prime}$.

Now suppose that $H^{\prime}\os H$. Let $N=\Phi(H)\times Z$. Then $N\ons G$.
Also $HN=H(\Phi(H)\times Z)=G$. Suppose for $K\cs G$ that $KN=G$. Now
$(KZ\cap H)\Phi(H)=H\cap (\Phi(H)KZ)=H\cap G=H$. So by standard
properties of the Frattini subgroup,
$KZ\cap H=H$. Thus by Proposition \ref{pgour},
$H/(K\cap H)\cong (KH\cap Z)/(K\cap Z)$. But the latter group is clearly
abelian and so $K\cap H\geqslant H^{\prime}$. As $H^{\prime}\os H$,
$K\cap H\os H$ and by the above isomorphism $K\cap Z\os KH\cap Z$. Suppose
that $K\neq H$. Then since $KN=G$, $K\geqslant H$ and so $KH\cap Z$ is
non-trivial. Hence $KH\cap Z\os Z$ and so $K\cap Z\os Z$ and $K\os G$. Thus
$H$ is solitary in $G$ by Lemma \ref{ldisub}(i).
\end{prf}

\begin{cor}
Let $G$ be an infinite finitely generated pro-$p$ group which is not
virtually $\mathbb Z_p$. Let $H\in\sub (G)^{\prime}$. If $H$ is central in $G$
then $H$ is not solitary in $G$.
\end{cor}

\begin{prf}
Suppose that $H$ is central in $G$ but that $H$ is solitary in $G$. Then
by Lemma \ref{ldisub}(vi) there exists $K\os G$ with $H\cns K$ and
$K/H\cong\mathbb Z_p$. Thus there exists $Z\cs G$ with $K=H\times Z$ and
$Z\cong\mathbb Z_p$. Now by Lemma \ref{ldisub}(ii) and Lemma \ref{lsdpzp},
$H^{\prime}\os H$. But $H$ is abelian and so $H$ is finite. Hence $G$ is
virtually $\mathbb Z_p$, a contradiction.
\end{prf}

\begin{cor}\label{cjizpi}
Let $H$ be a non-abelian just-infinite pro-$p$ group, $Z\cong\mathbb Z_p$
and $G=H\times Z$. Then $G$ has precisely countably infinitely many solitary
subgroups.
\end{cor}

\begin{prf}
By Lemma \ref{lctbdx}(iii), $G$ has at most countably infinitely many
solitary subgroups. Let $N\ons H$. Then $N^{\prime}\ons H$ since $H$
is non-abelian and just-infinite. By Lemma \ref{lsdpzp}, $N$ is solitary in
$N\times Z$. But $N\times Z\os G$ and so by Lemma \ref{ldisub}(ii), $N$ is
solitary in $G$. Since $H$ has precisely countably infinitely many open
normal subgroups, we are done.
\end{prf}

By taking $H$ to be a non-abelian soluble just-infinite pro-$p$ group in
Corollary \ref{cjizpi} we obtain examples of poly-procyclic pro-$p$ groups
with countably infinitely many solitary subgroups. Even poly-$\mathbb Z_p$
groups $G$ with $\h (G)>1$ may not have $\sub (G)$ homeomorphic to
Pe\l czy\'nski space.

\begin{propn}\label{pirim}
Let $G$ be a pro-$p$ group. Suppose $G$ has a closed normal subgroup $H$
with $H\cong\mathbb Z_p^n$, for $n$ an integer $>1$, and
$G/H\cong\mathbb Z_p$. Suppose further that $H$ is a rationally irreducible
$L$-module for every $L\os G$, with $H\leqslant L$. (This means that if
$K\cns L$ with $K\leqslant H$ then either $K=1$ or $K\ons H$.) Then $H$ is
solitary in $G$.
\end{propn}

\begin{prf}
$\Phi (H)\os H$ so (Proposition~1.2(v), \cite{Dix}) there exists $M\os G$ such
that $M\cap H=\Phi (H)$. Let $N=M_G$. So $N\ons G$. Suppose $HN=KN$ where
$K\in\sub (G)^{\prime}$. By Lemma \ref{ldisub}(i) it suffices to show that
$K=H$. Firstly we claim that $H\cap K$ is non-trivial. So suppose for a
contradiction that $H\cap K=1$. If $K\leqslant H$ then $K=1$, and so
$H\leqslant N$. Thus $H=H\cap N\leqslant H\cap M=\Phi (H)$, a
contradiction. Hence $K\not\leqslant H$. Thus $HK/H\cong\mathbb Z_p$ and
so $K\cong\mathbb Z_p$. Consequently $HN/N=KN/N$ is cyclic and so
$H/(H\cap N)$ is cyclic. Hence $H/\Phi (H)=H/(M\cap H)$ being a homomorphic
image of $H/(H\cap N)$ is cyclic, a contradiction as $n>1$. Thus
$H\cap K\neq 1$.

Now suppose for a contradiction that $K\nleqslant H$. As $H$ is abelian,
$H\cap K\cns H$ and as $H\cns G$, $H\cap K\cns K$. Thus
$HK\leqslant\norm_G(H\cap K)$. But as $K\nleqslant H$, $HK\os G$ and so
$\norm_G(H\cap K)\os G$. By hypothesis $H$ is a rationally irreducible
$(\norm_G(H\cap K))$-module. So as $H\cap K\neq 1$, $H\cap K\ons H$. Now
there is a non-negative integer $r$ such that $H\cap K=H^{p^r}$. So
$K\geqslant H^{p^r}$. But $K\nos G$ so $K\leqslant H$, a contradiction.
Hence $K\leqslant H$. Now
$H=H\cap HN=H\cap KN=(H\cap N)K\leqslant (H\cap M)K=\Phi (H)K$. So
$H=\Phi (H)K$ and $K=H$ by standard properties of the Frattini
subgroup, as required. 
\end{prf}

We now give an example of a group satisfying the hypotheses of Proposition
\ref{pirim}, for simplicity in the case where $p=3$ and $n=2$.

\begin{example}
There are poly-$\mathbb Z_p$
groups $G$ with scattered height $2$ and infinitely many solitary subgroups..
\end{example}

\begin{prf}
Every extension by $\mathbb Z_3$ splits (since $\mathbb Z_3$ is a free
pro-$3$ group). So
we must specify a group of the form $G=H\rtimes T$ where
$H\cong\mathbb Z_3\times\mathbb Z_3$ and
$T=\overline{\langle t\rangle}\cong\mathbb Z_3$. Consequently we must
specify a continuous homomorphism from $T$ to $\Aut (H)$, and since $T$
is a free pro-$3$ group and $\Aut (H)\cong\gl{2}{\mathbb Z_3}$, it
suffices to specify a matrix in a pro-$3$ subgroup of
$\gl{2}{\mathbb Z_3}$ to which $t$ is mapped. Let $\alpha=1+\sqrt 5$.
Consider the field extension $\mathbb Q_3(\alpha):\mathbb Q_3$ where
$\mathbb Q_3$ is the field of $3$-adic numbers; that is, the field of
fractions of $\mathbb Z_3$. We consider $\mathbb Q_3(\alpha)$ as a
vector space over $\mathbb Q_3$. Multiplication (on the right) by
$\alpha$ is an vector space endomorphism of $\mathbb Q_3(\alpha)$. Let
$A=\left(\begin{array}{cc}0 & 1\\4&2\end{array}\right)$. Then $A$ is
the matrix corresponding to this endomorphism with respect to the basis
$\{1,\alpha\}$. Let
$K_1=\{g\in\gl{2}{\mathbb Z_3}\mid g-1\in 3M_2(\mathbb Z_3)\}$. This is
the \df{first congruence subgroup} of $\gl{2}{\mathbb Z_3}$, and is an
open normal pro-$3$ subgroup of $\gl{2}{\mathbb Z_3}$; see, for example,
\S8.5 of \cite{WilPrf}. We choose a power $A_1$ of $A$ lying in $K_1$
and define $G$ as above by mapping $t$ to $A_1$.

Now for $n$ a positive integer, $\alpha^n=a_n+b_n\alpha$ for some
$a_n,b_n\in\mathbb Q_3$. Clearly if $n=1$ or $2$, then $b_n>0$. It is
straightforward to show that if $n>2$ then $b_n=2b_{n-1}+4b_{n-2}$. So by
induction, $b_n>0$ for every positive integer $n$. Thus
$\alpha^n\notin\mathbb Q_3$ for every $n$. It also follows that
$(1-\sqrt 5)^n\notin\mathbb Q_3$ for every $n$. The eigenvalues of $A$
are $1\pm\sqrt 5$. For each $n$, the eigenvalues of $A^n$ are the $n$th
power of the eigenvalues of $A$, namely $(1\pm\sqrt 5)^n$. Hence $A^n$ has
no eigenvalues in $\mathbb Q_3$ for each positive integer $n$. Consequently
$A_1^n$ also has no eigenvalues in $\mathbb Q_3$ for each positive integer
$n$.

We want to show that every open subgroup of $T$ acts rationally irreducibly
on $H$. Suppose for a contradiction that there exists $L\os T$ and $K\cns H$
with $K\neq 1$ and $K\nos H$. Now $L=T^{3^r}$ for some non-negative integer
$r$. $L$ acts on $H$ by mapping $t^{3^r}$ to $A_1^{3^r}$. Then
$K\otimes_{\mathbb Z_3}\mathbb Q_3$ is a $1$-dimensional
$\mathbb Q_3$-subspace of $H\otimes_{\mathbb Z_3}\mathbb Q_3$. So $t^{3^r}$
acts as scalar multiplication on $K\otimes_{\mathbb Z_3}\mathbb Q_3$ and
thus the matrix $A_1^{3^r}$ has an eigenvalue in $\mathbb Q_3$, a
contradiction.
\end{prf}

\paragraph{Scattered Height Equal to $2$, and Finitely Many Isolated
  Subgroups -- $G$   Virtually $\mathbb{Z}_p$} 

In Corollary \ref{cjizpi} we gave examples of pro-$p$ groups with
infinitely many solitary subgroups. There also exist profinite groups $G$
with $\sub (G)$ not homeomorphic to Pe\l czy\'nski space but with only
finitely many solitary subgroups. For example if $G$ is a pro-$p$ group
with $\sub (G)$ countably infinite, then (see \cite{GS1})
$\sub(G)^{\prime}$ is a finite discrete space.

\begin{lemma}\label{lisvzp}
Let $G$ be a profinite group. Then $G$ is virtually $\mathbb Z_p$ for some
prime $p$ if and only if $G$ has a finite solitary subgroup. Also if $G$
is virtually $\mathbb Z_p$ for some prime $p$ then every
$H\in\sub (G)^{\prime}$ is finite.
\end{lemma}

\begin{prf}
By Lemma \ref{ldisub}(vii) if $G$ is virtually $\mathbb Z_p$ for some prime
$p$ then $1$ is solitary in $G$. Conversely suppose that $G$ has a finite
and solitary subgroup $H$. By Lemma \ref{ldisub}(vi) there exists $K\os G$ with
$H\cns K$ and $K/H\cong\mathbb Z_p$ for some prime $p$. So there exists
$Z\cs G$ with $HZ=K$, $H\cap K=1$ and $Z\cong\mathbb Z_p$. But
$|K:Z|=|HZ:Z|=|H|$ which is finite by hypothesis. Thus $Z\os G$ and $G$ is
virtually $\mathbb Z_p$.

Now suppose that $G$ is virtually $\mathbb Z_p$ for some prime $p$. Let
$Z\ons G$ with $Z\cong\mathbb Z_p$. Let $H\in\sub (G)^{\prime}$. If $H\cap Z$
is non-trivial then $H\cap Z\os Z$ and so $H\os G$, a contradiction. So
$H\cap Z=1$. Hence $|H|=|H:H\cap Z|=|HZ:Z|$ which is finite.
\end{prf}

\begin{propn}\label{pprvzpc}
Let $G$ be a profinite group with an open normal subgroup $Z$, topologically
isomorphic to $\mathbb Z_p$ for some prime $p$. Let $H\in\sub (G)^{\prime}$.
Then $H$ is solitary in $G$ if and only if $H\leqslant \cen_G(Z)$. Moreover
there are only finitely many such subgroups.
\end{propn}

\begin{prf}
Clearly $Z\leqslant\zg (\cen_G(Z))$. So $\cen_G(Z)$ has an open central
subgroup topologically isomorphic to $\mathbb Z_p$. Hence by
Theorem~4.2 
of \cite{GS1}, $\sub (\cen_G(Z))^{\prime}$ is finite. Now if
$H\leqslant \cen_G(Z)$ then $H\in\sub (\cen_G(Z))^{\prime}$ and so $H$ must
be solitary in $\cen_G(Z)$. Also $\cen_G(Z)\ons G$. So by Lemma
\ref{ldisub}(ii) $H$ is solitary in $G$.

Now suppose that $H$ is solitary in $G$. We know by Lemma \ref{lisvzp} that
$H$ is finite and by the proof of Lemma \ref{lisvzp} that $H\cap Z=1$. For
each $h,k\in H$ the maps $g\mapsto h^g$ and $g\mapsto k$ are continuous
maps from $\norm_G(H)$ to $G$. Thus for each $h,k\in H$,
$\{g\in\norm_G(H)\mid h^g=k\}$ is closed in $\norm_G(H)$. Hence for each
$h\in H$ we can write $\norm_G(H)=\bigcup_{k\in H}\{g\in\norm_G(H)\mid h^g=k\}$
a finite union of closed sets. So by Baire's category theorem for each $h\in H$
there exists $k_h\in H$ such that $\{g\in\norm_G(H)\mid h^g=k_h\}$ has non-empty
interior in $\norm_G(H)$. But by Lemma \ref{ldisub}(vi) $\norm_G(H)\os G$ and
so this set has non-empty interior in $G$. Consequently (as the cosets
of open normal subgroups form a base for the topology on $G$) 
for each $h\in H$ there exists $k_h\in H$, $g_h\in G$
and $N_h\ons G$ such that $\{g\in\norm_G(H)\mid h^g=k_g\}\supseteq N_hg_h$. Let
$h\in H$ and $n\in N_h$. Then $h^{ng_h}=k_h$ and so $h^n=k_h^{g_h^{-1}}$ which
is independent of $n$. Consequently, $h^n=h$ for every $n\in N_h$.

Now let $N=Z\cap\bigcap_{h\in H}N_h$. Then since $H$ is finite, $N\ons G$.
Also as $H\cap Z=1$, $N\cap Z=1$. Let $h\in H$ and $n\in N$. Then as
$n\in N_h$, $[h,n]=1$. Thus $H$ acts trivially on $N$ by conjugation. Also
$H$ acts on $Z$. $Z$ has a natural ring structure, and 
 $\Aut (Z)\cong U(Z)$, the group of units of $Z$. Let $h\in H$.
Then there exists $u\in U(Z)$ such that $u\cdot z=z^h$ for every $z\in Z$
where $\cdot$ is the ring multiplication in $Z$. But $H$ acts trivially
on $N$ and so $u\cdot n=n$ for every $n\in N$. If we choose a non-trivial
(i.e. non-zero in the ring) $n\in N$ then we may cancel it since $N$ is an
integral domain and conclude that $u=1$. Thus $z^h=z$ for every $h\in H$ and
$z\in Z$. Hence $H\leqslant\cen_G(Z)$ as required.
\end{prf}

\begin{cor}\label{cprvztop}
Let $G$ be a profinite group which is virtually $\mathbb Z_p$ for some
prime $p$. 
Suppose that $\sub (G)$ is uncountable. Then $\sub (G)$ is homeomorphic to
$P \oplus (\omega n+1)$ where
$n=|\sub (\cen_G(Z))^{\prime}|$ (and $P$ is the Pe\l czy\'nski space).
\end{cor}

\begin{prf}
The result follows immediately from Proposition \ref{pprvzpc} and Proposition
\ref{pchpelconv}.
\end{prf}

\begin{cor}
Let $G$ be a just-infinite non-abelian pro-$p$ group which is virtually
$\mathbb Z_p$. Then $\sub (G)$ is homeomorphic to
$P\oplus (\omega +1)$.
\end{cor}

\begin{prf}
Suppose $G$ has non-trivial centre. Then by Proposition~3.6(i) 
 of \cite{GS1},
$G^{\prime}$ is finite. Thus as $G$ is just-infinite, $G^{\prime}$ is trivial;
that is, $G$ is abelian, a contradiction. Hence $G$ has trivial centre and so by
Theorem~3.7 
 of \cite{GS1}, $\sub (G)$ is uncountable. Let $Z$
be a maximal normal 
abelian subgroup of $G$. Then $Z\cns G$ and since $G$ is just-infinite,
$Z\ons G$. Since $G$ is virtually $\mathbb{Z}_p$ it has solitary
subgroups, and so definitely isolated subgroups. So  $G$ 
 and its open subgroup $Z$ are finitely generated. Also as $G$ is
just-infinite, $\T (Z)=1$. Thus as $G$ 
is virtually $\mathbb Z_p$, $Z$ must be topologically isomorphic to
$\mathbb Z_p$. Now $G/Z$ is nilpotent and so by a standard argument (see for
example the proof of 5.2.3 of \cite{Rob}), $Z=\cen_G(Z)$. The result now
follows from Corollary \ref{cprvztop}.
\end{prf}

\paragraph{Scattered Height Equal to $2$, and Finitely Many Solitary
  Subgroups -- $G$   Not Virtually $\mathbb{Z}_p$} 

It is natural to ask whether an infinite finitely generated pro-$p$ group
which is not virtually $\mathbb Z_p$ can have solitary subgroups but have only
finitely many solitary subgroups. It is clear from the proof of Corollary
\ref{cjizpi} that it might not be possible to find such examples if it is the
case that if an infinite finitely generated pro-$p$ group has an open
commutator subgroup then it has a proper open subgroup with an open commutator
subgroup. This though does not happen.

\begin{example}[Wilson]\label{exwil}
Let $G$ be the pro-$2$ group with presentation
$\langle x_1,x_2\mid (x_1^2)^{x_2}=x_1^{-2},(x_2^2)^{x_1}=x_2^{-2},
((x_1x_2)^2)^{x_1}=(x_1x_2)^{-2} \rangle$. Then $G$ is a poly-procyclic pro-$2$
group. $G^{\prime}\os G$ but for every proper open subgroup $H$ of $G$,
$H^{\prime}\nos H$.
\end{example}

\begin{prf}
We sketch a proof that $G$ has these properties. For the definition of pro-$2$
presentations, see \S 12.1 of \cite{WilPrf}. In particular $G$ is the pro-$2$
completion of the corresponding abstract group with this presentation.

Let $a_1=x_1^2, a_2=x_2^2, a_3=(x_1x_2)^2$, and
$A=\overline{\langle a_1,a_2,a_3\rangle}$. Since
$a_3=a_3^{x_1x_2}=(a_3^{-1})^{x_2}$, $a_3^{x_2}=a_3^{-1}$. Thus $A\cns G$. Also
it is easy to check that $[a_1,a_2]=[a_1,a_3]=[a_2,a_3]=1$. So $A$ is abelian.
Now each $a_i$ is of infinite order. This can be seen, for example, by embedding
$G$ into the direct product of three copies of the pro-$2$ analogue of the
infinite dihedral group. Consequently the $a_i$'s
are distinct and $A\cong\mathbb Z_2\times\mathbb Z_2\times\mathbb Z_2$. Clearly
$G/A$ is a Klein $4$-group. So $G$ is a poly-procyclic pro-$2$ group
of Hirsch length $3$.

Now $[x_2,a_1]=a_1^2$, $[x_1,a_2]=a_2^2$, and $[x_1,a_3]=a_3^2$. So
$\overline{\langle a_1^2,a_2^2,a_3^2\rangle}\leqslant G^{\prime}$. Thus
$G^{\prime}\os G$. Since $G^2=A$, and $G^{\prime} \cs A$ ($G/A$
abelian), 
we see that
$\Phi (G)=A$ (for any pro-$p$ 
group $G$ we have $\Phi(G) = \overline{G^{\prime} G^p}$, see Proposition~2.5.2
\cite{WilPrf}).
Thus the maximal open subgroups of $G$ are $\langle x_1,A\rangle$,
$\langle x_2,A\rangle$, and $\langle x_1x_2,A\rangle$. It is now easy to check
that $\langle x_1,A\rangle^{\prime}=\overline{\langle a_2^2,a_3^2\rangle}$,
$\langle x_2,A\rangle^{\prime}=\overline{\langle a_1^2,a_3^2\rangle}$, and
$\langle x_1x_2,A\rangle^{\prime}=\overline{\langle a_1^2,a_2^2\rangle}$. Now
each of these groups is topologically isomorphic to
$\mathbb Z_2\times\mathbb Z_2$; in particular they each have Hirsch length $2$.
Let $H$ be a proper open subgroup of $G$. Then $H$ is contained in a maximal
open subgroup $M$ of $G$. Now $H^{\prime}\leqslant M^{\prime}$ and so
$\h (H^{\prime})\leqslant 2 < 3=\h (G)$. Hence $H^{\prime}\nos H$.
\end{prf}

\begin{propn}\label{pdiprwzp}
Let $W$ be an infinite pro-$p$ group of finite rank, $Z\cong\mathbb Z_p$ and
$G=W\times Z$. Let $H\in\sub (G)^{\prime}$.
\begin{enumerate}
\item[(i)] Suppose that $H$ is solitary in $G$. Then either $H\os W$ or
$H\cap Z\neq 1$. If $H\cap Z\neq 1$ then $H^{\prime}\nos H$. In either case
$(H\cap W)\times (H\cap Z)\os H\os (HZ\cap W)\times (HW\cap Z)$.
\item[(ii)] Suppose $H$ factorises, i.e. $H=(H\cap W)\times (H\cap Z)$. Then
$H$ is solitary in $G$ if and only if $H\os W$ and $H^{\prime}\os H$.
\end{enumerate}
\end{propn}

\begin{prf}
For (i), suppose that $H$ is solitary in $G$. Then by Lemma
\ref{ldisub}(vi), for some open subgroup $K$ of $G$, $H \ons K$ and
$K/N \cong \mathbb{Z}_p$. Hence $\dim (H)= \dim (K) - \dim
(\mathbb{Z}_p) = \dim (G)-1 = \dim (W)$.

Firstly suppose that $H\cs W$. Since $\dim (H)= \dim(W)$,  $H\os W$. Clearly
$(H\cap W)\times (H\cap Z)\os H\os (HZ\cap W)\times (HW\cap Z)$. Now suppose
that $H\not\leqslant W$, that is that $HW\cap Z$ is non-trivial. Suppose
for a contradiction that $H\cap Z$ is trivial. Then
$(HW\cap Z)/(H\cap Z)\cong\mathbb Z_p$. Thus by Proposition \ref{pgour},
$$
\frac{(HZ\cap W)\times (HW\cap Z)}{(H\cap W)\times (H\cap Z)}
\cong\mathbb Z_p\times\mathbb Z_p.
$$
If $(HZ\cap W)\times (HW\cap Z)\os G$ then by Lemma \ref{ldisub}(ii,iii),
$$
\frac{H}{(H\cap W)\times (H\cap Z)}\hbox{ is isolated in }
\sub\left(\frac{(HZ\cap W)\times (HW\cap Z)}{(H\cap W)\times (H\cap Z)}
\right)^{\prime},
$$
contradicting Examples \ref{ezpnfree}(i).
Thus $(HZ\cap W)\times (HW\cap Z)$ is not open $G$, so
$\dim ((HZ\cap W)\times (HW\cap Z)) \le \dim (W) = \dim (G)-1$. 
From the above isomorphism
$\dim (H\cap W) \le \dim (W)-2$. But by Proposition \ref{pgour}
$H/((H\cap W)\times (H\cap Z))\cong\mathbb Z_p$. So
$\dim (H\cap W)=\dim(H) - \dim (\mathbb{Z}_p) = \dim (W)-1$, a
contradiction. Thus $H\cap Z$ is non-trivial.

Now from Proposition \ref{pgour}
$(H\cap W)\times (H\cap Z)\os H\os (HZ\cap W)\times (HW\cap Z)$. Note
that since $HW\cap Z\cong\mathbb Z_p$,
$HZ\cap W\nos (HZ\cap W)\times (HW\cap Z)$. Suppose that
$H^{\prime}\os H$. Then $H^{\prime}\os (HZ\cap W)\times (HW\cap Z)$. Hence
$((HZ\cap W)\times (HW\cap Z))^{\prime}\os (HZ\cap W)\times (HW\cap Z)$.
But $((HZ\cap W)\times (HW\cap Z))^{\prime}=(HZ\cap W)^{\prime}$. So
$HZ\cap W\os (HZ\cap W)\times (HW\cap Z)$, a contradiction. Hence
$H^{\prime}\nos H$.

For (ii), suppose that $H$ factorises and that $H$ is solitary in $G$.
Suppose for a contradiction that $H\nos W$. By Lemma \ref{ldisub}(vi)
there exists $K\os G$ with $H\cns K$ and $K/H\cong\mathbb Z_p$. Now
$$
\frac{K\cap W}{H\cap W}\cong\frac{H(K\cap W)}{H}\leqslant\frac{K}{H}.
$$
If $H=H(K\cap W)$ then $K\cap W\leqslant H$ and so $H\cap W=K\cap W$. But
$\dim (H\cap W)=\dim (W)-1$ where as $\dim (K\cap W)=\dim (W)$, a
contradiction. Thus $H\neq H(K\cap W)$ and so
$(K\cap W)/(H\cap W)\cong\mathbb Z_p$. Now $(K\cap W)\times Z\os G$ and
$H=(H\cap W)\times (H\cap Z)\leqslant (K\cap W)\times Z$. So by Lemma
\ref{ldisub}(ii,iii),
$$
\frac{H}{H\cap W}\hbox{ is isolated in }\sub\left(\frac{K\cap W}{H\cap W}
\times Z\right)^{\prime}.
$$
But by Examples \ref{ezpnfree}(i) this is a contradiction since
${\displaystyle\frac{K\cap W}{H\cap W}\times Z\cong\mathbb Z_p\times
\mathbb Z_p}$. Thus $H\os W$.

Now by Lemma \ref{ldisub}(ii) $H$ is solitary in $H\times Z$. So by Lemma
\ref{lsdpzp}, $H^{\prime}\os H$.

Conversely suppose that $H$ factorises, $H\os W$ and $H^{\prime}\os H$. By
Lemma \ref{lsdpzp} $H$ is solitary in $H\times Z$ and so by Lemma
\ref{ldisub}(ii), $H$ is solitary in $G$.
\end{prf}

\begin{cor}
Let $W$ be the group of Example \ref{exwil}, $Z\cong\mathbb Z_2$, and
$G=W\times Z$. Let $H\in\sub (G)^{\prime}$. If $H$ is solitary in $G$ then
either $H=W$ or $H\cap Z\neq 1$ and $H^{\prime}\nos H$. If $H$ factorises
then $H$ is solitary in $G$ if and only if $H=W$.
\end{cor}

\begin{prf}
Suppose that $H$ is solitary in $G$. If $H\cs W$ then by Proposition
\ref{pdiprwzp}(i), $H\os W$. So by Lemma \ref{ldisub}(ii) and Lemma
\ref{lsdpzp} $H^{\prime}\os H$. But then by the properties of $W$, $H=W$.
If $H\nleqslant W$ then by Proposition \ref{pdiprwzp}(i) $H\cap Z\neq 1$
and $H^{\prime}\nos H$.

Now suppose that $H$ factorises and is solitary in $G$. Then by Proposition
\ref{pdiprwzp}(ii) and the above, clearly $H=W$. Conversely if $H=W$ then
again by Proposition \ref{pdiprwzp}(ii) $H$ is solitary in $G$ since
$W^{\prime}\os W$.
\end{prf}

\begin{prob}\label{pbnvzpfis}
Can an infinite finitely generated pro-$p$ group which is not virtually
$\mathbb Z_p$, have solitary subgroups, but have only finitely many
solitary subgroups?

Does the group given in Example~\ref{exwil} above,  provide a
positive answer to this problem?
\end{prob}

\section{General Profinite Groups Attain All Permitted Scattered
  Heights}\label{gen_case} 

We have seen above that pro-$p$ groups can have subgroup spaces of all
heights permitted by Theorem~\ref{bd_sc_ht} -- namely scattered
height $1$ or $2$.  Now we show that all scattered heights permitted
by Theorem~\ref{bd_sc_ht} for a general profinite group can also be
attained.

\begin{thm}
For each $n \le k+1$, there is a profinite group $G(n,k)$ satisfying:

the number of primes $p$ such that $p^\infty$ divides $o(G(n,k))$ equals
$k$ and $\hht(\sub(G)) = n$.
\end{thm} 
\begin{prf}
Fix $k$. We define appropriate profinite groups $G(n,k)$ for $n=1$,
$1<n \le k$ and $n=k+1$. We choose $G(n,k)$  equal to $\Dr_{i=1}^k
G_{p_i}$ where the $p_i$s are distinct primes. In this case $\sub
(G(n,k)) \cong \prod_{i=1}^k \sub ( G_{p_i})$.

\medskip

\noindent {\sl Case $n=1$:} \ Let $G_{p_i}= \mathbb{Z}_{p_i}^2$. Then
$\sub(G_{p_i})=P$ (the Pe\l czy\'nski space), and $\sub (G(n,k) \cong
P^k \cong P$, which has scattered height $1(=n)$.

\medskip 

\noindent {\sl Case $1<n \le k$:} \ Let $G_{p_i}= \mathbb{Z}_{p_i}$
for $i=1, \ldots , n$ and $G_{p_i}= \mathbb{Z}_{P-I}^2$ for $i=n+1,
\ldots , k+1$. Then $\sub(G(n,k))\cong (\omega+1)^n \times P^{k+1-n}$,
which has scattered height $n$. 

\medskip

\noindent {\sl Case $n= k+1$:} \ Let $G_{p_i}= \mathbb{Z}_{p_i}$
for $i=1, \ldots , k$ and $G_{p_{k+1}}$ be any pro-$p_{k+1}$ group
with uncountable subgroup space $X$ with scattered height $2$ -- for
example a non-abelian just-infinite pro-$p_{k+1}$ direct producted
with $\mathbb{Z}_{p_{k+1}}$.  Then
$\sub(G(n,k)) \cong (\omega+1)^k \times X$, 
which has scattered height $k+1(=n)$.
\end{prf}

\bibliographystyle{plain}

\bibliography{sgbib}

\end{document}